\documentclass[12pt]{amsart}

\usepackage{amsmath,amssymb,amsthm,amscd,a4wide,upref}

\usepackage[enableskew]{youngtab}

\usepackage[colorlinks=true, linkcolor=blue, citecolor=magenta, urlcolor=blue, backref=page]{hyperref} 

\hypersetup{colorlinks,citecolor=blue,filecolor=black,linkcolor=blue,urlcolor=blue}
\usepackage{enumerate}
\usepackage{mathtools}
\usepackage{enumitem}

\usepackage{lmodern}     
\usepackage[T1]{fontenc}

\newcommand{\End}{{\rm{End}\ts}}

\newcommand{\diag}{ {\rm diag}}

\newcommand{\wh}{\widehat}
\newcommand{\ot}{\otimes}

\newcommand{\De}{\Delta}

\newcommand{\si}{\sigma}

\newcommand{\de}{\delta}

\newcommand{\om}{\omega}

\newcommand{\ts}{\,}
\newcommand{\tss}{\hspace{1pt}}

\newcommand{\qin}{q^{-1}}

\newcommand{\U}{ {\rm U}}

\newcommand{\CC}{\mathbb{C}\tss}

\newcommand{\Y}{ {\rm Y}}

\newcommand{\ol}{\overline}

\newcommand{\Hc}{\mathcal{H}}

\newcommand{\gl}{\mathfrak{gl}}

\newcommand{\z}{\mathfrak{z}}

\newcommand{\str}{ {\rm str}}

\newcommand{\bi}{\bar{\imath}}
\newcommand{\bj}{\bar{\jmath}}
\newcommand{\bk}{\bar{k}}

\newcommand{\Sym}{\mathfrak S}

\newcommand{\Fand}{\qquad\text{and}\qquad}

\newtheorem{thm}{Theorem}[section]
\newtheorem{lem}[thm]{Lemma}
\newtheorem{prop}[thm]{Proposition}
\newtheorem{cor}[thm]{Corollary}
\newtheorem{conj}[thm]{Conjecture}

\theoremstyle{definition}
\newtheorem{defin}[thm]{Definition}

\theoremstyle{remark}
\newtheorem{remark}[thm]{Remark}
\newtheorem{example}[thm]{Example}

\newcommand{\bth}{\begin{thm}}
\renewcommand{\eth}{\end{thm}}
\newcommand{\bpr}{\begin{prop}}
\newcommand{\epr}{\end{prop}}
\newcommand{\ble}{\begin{lem}}
\newcommand{\ele}{\end{lem}}
\newcommand{\bco}{\begin{cor}}
\newcommand{\eco}{\end{cor}}
\newcommand{\bde}{\begin{defin}}
\newcommand{\ede}{\end{defin}}
\newcommand{\bex}{\begin{example}}
\newcommand{\eex}{\end{example}}
\newcommand{\bre}{\begin{remark}}
\newcommand{\ere}{\end{remark}}
\newcommand{\bcj}{\begin{conj}}
\newcommand{\ecj}{\end{conj}}

\newcommand{\bal}{\begin{aligned}}
\newcommand{\eal}{\end{aligned}}
\newcommand{\beq}{\begin{equation}}
\newcommand{\eeq}{\end{equation}}
\newcommand{\ben}{\begin{equation*}}
\newcommand{\een}{\end{equation*}}

\newcommand{\bpf}{\begin{proof}}
\newcommand{\epf}{\end{proof}}

\newcommand{\BC}{\mathbb{C}}            
\newcommand{\BZ}{\mathbb{Z}}            

\def\beql#1{\begin{equation}\label{#1}}

\begin{document}
\title{Quantum Berezinian for quantum affine superalgebra $\mathrm{U}_q(\widehat{\gl}_{M|N})$}
\author{Naihuan Jing}
\address{Department of Mathematics, North Carolina State University, Raleigh, NC 27695, USA}
\email{jing@ncsu.edu}
\author{Zheng Li}
\address{School of Mathematics and Statistics, Central China Normal University, Wuhan, Hubei 430079, China}
\email{lz1994@ccnu.edu.cn}
\author{Jian Zhang}
\address{School of Mathematics and Statistics, and Hubei Key Lab-Math. Sci., Central China Normal University, Wuhan, Hubei 430079, China}
\email{jzhang@ccnu.edu.cn}

\thanks{{\scriptsize
\hskip -0.6 true cm MSC (2010): Primary: 17B37; Secondary: 81R50.
}}
\maketitle

\begin{abstract}
We introduce the quantum Berezinian for the quantum affine superalgebra $\mathrm{U}_q(\widehat{\mathfrak{gl}}_{M|N})$ and
show that the coefficients of the quantum Berezinian belong to the center of $\mathrm{U}_q(\widehat{\gl}_{M|N})$.
We also construct another family of central elements which can be expressed in the quantum Berezinian by a Liouville-type theorem.
Moreover, we prove analogues of the Jacobi identities, the Schur complementary theorem, the Sylvester theorem and the MacMahon Master theorem for the generator matrices of $\mathrm{U}_q(\widehat{\gl}_{M|N})$.
\end{abstract}

\section{Introduction}
\label{sec:int}
\setcounter{equation}{0}

Let $\gl_{M|N}$ be the general linear Lie superalgebra over $\mathbb C$,
there is a natural homomorphism from the Yangian $\Y(\gl_{M|N})$ to the enveloping superalgebra
$\mathrm{U}(\gl_{M|N})$. Nazarov \cite{Na} showed that the coefficients of the Berezinian of the matrix $(t_{ij}(z+(M-N-j))_{1\leq i,j\leq M}\oplus (t^*_{ij}(z-(M+N-j+1))_{M+1\leq i,j\leq M+N}$ belong to the center of $\Y(\gl_{M|N})$, therefore naturally project to the center of $\mathrm{U}(\gl_{M|N})$.
This gives rise to the Capelli identity in $\Y(\gl_{M|N})$ and $\mathrm U(\gl_{M|N})$ \cite{Na, MR, Na2}. Gow \cite{G}
showed that the center of $\Y(\gl_{M|N})$ is generated by the coefficients of the Berezinian.
Further properties of the Berezinian and $\Y(\mathfrak{gl}_{M|N})$ were investigated in \cite{P, LM, T, CH}.
Nazarov also constructed another family of central elements by deriving the Liouville-type formulas, and this family is also related to
the Berezinian. Recently Bagnoli and Kozic also constructed central elements via the Berezinian for the double Yangian $DY(\mathfrak{gl}_{M|N})$ \cite{BK}.

The quantum superalgebra $\mathrm{U}_q(\widehat{\gl}_{M|N})$ \cite{FHS, GZ, Y} is a deformation of the enveloping algebra $\mathrm U(\widehat{\gl}_{M|N})$ of the affine Lie superalgebra
$\widehat{\gl}_{M|N}$, its specialization at the trivial central charge is also called the quantum loop superalgebra \cite{Z, KwL}.
The associative algebra $\mathrm{U}_q(\widehat{\gl}_{M|N})$ is
generated by the quantum root vectors ${l_{ij}^{\pm}}^{{(r)}}$ subject to the relations \eqref{RLL-0}, \eqref{RLL-1}, \eqref{RLL} and \eqref{RLL cros}. The
generators ${l_{ij}^{\pm}}^{{(r)}}$ can be also defined as the coefficients of the generator series $l_{ij}^{\pm}(z)=\sum_{r=0}^{\infty}{l_{ij}^{\pm}}^{(r)}z^{\mp r}$ subject
the RLL relations governed by the quantum Yang-Baxter equation \eqref{QYBE}. Quantum affine superalgebras have found applications in physical models and
also appeared in the AdS/CFT correspondence \cite{LRT}. On the other hand,
the Yangian $\mathrm Y(\gl_{M|N})$ \cite{MR, GTL, LWZ} and double Yangian $DY(\mathfrak{gl}_{M|N})$ \cite{Zh, BK}
can be viewed as some limits of the quantum affine superalgebra
$\mathrm{U}_q(\widehat{\gl}_{M|N})$, therefore one may expect that important structures may be lifted to the quantum affine superalgebra $\mathrm{U}_q(\widehat{\gl}_{M|N})$.

In this paper, we will construct two families of central elements
for $\mathrm{U}_q(\widehat{\gl}_{M|N})$ ($M\neq N$) which are quantum analog of the central elements for the super Yangian given by Nazarov.
Our first result is the following description of the quantum Berezinian:
\beql{ber}
\bal
B_q
&(L^{\pm}(z) )
=\sum_{\sigma\in S_M}(-q)^{-l(\sigma)}
L^{\pm}(z)_{ {\sigma(1)}, 1}\cdots L^{\pm}(zq^{2M-2})_{ {\sigma (M)},M}\\
&\times \sum_{\tau\in S_N} (-q)^{-l(\tau)}  \big({L^{\pm}(zq^{2M-2}) }^{-1}\big)_{M+1,M+\tau(1)}
\cdots \big({L^{\pm}(zq^{2M-2N})}^{-1}\big)_{M+N,{M+\tau(N)}}\\
&=b^{\pm}_0+b^{\pm}_1z^{\mp 1}+b^{\pm}_2z^{\mp 2}+\cdots
\eal
\eeq

Then we show that the coefficients $b^{\pm}_i (i\in\mathbb Z_+)$ of the quantum Berezinian belong to the center of  $\U_q(\wh{\gl}_{M|N})$.
  For simplicity we have adopted the standard parity in defining $\mathrm{U}_q(\widehat{\gl}_{M|N})$ and its quantum Berezinian. We remark that
one can define both for any parity sequence $s$ of $M$ zeros and $N$ ones just as in the super Yangian case (cf. \cite{P, T, Y}), however, some formulas will be more involved. 

We also show that the Liouville theorem holds in the quantum affine superalgebra $\U_q(\wh{\gl}_{M|N})$
which provides another family of central elements. Let $D$ be the following diagonal element
\begin{equation}
D=\diag\big[q^{2},q^{4},\ldots,q^{2M},q^{2M},\ldots,q^{2M-2N+2}\big].
\end{equation}

\bpr\label{prop:equiz}
There exist a series $\z^+(z)$ in $z^{-1}$ and a series $\z^-(z)$ in $z$ with coefficients
in the quantum affine superalgebra $\mathrm{U}_q(\widehat{\gl}_{M|N})$ such that
\begin{align}\label{lltr}
	L^{\pm}(zq^{2N-2M})^{st}D(L^{\pm}(z)^{-1})^{st}&=\z^{\pm}(z)D \\
	\intertext{and}
	(L^{\pm}(z)^{-1})^{st}D^{-1} L^{\pm}(zq^{2N-2M})^{st}&=\z^{\pm}(z)D^{-1}.
	\label{secllrt}
\end{align}
Here the superscript $st$ means taking supertranspose. Moreover, the coefficients of the series $\z^{\pm}(u)$ belong to the center of the algebra
$\U_q(\wh{\gl}_{M|N})$.
\epr

In fact, the two families of central elements are related with the help of the Berezinian as we have the following Liouville-type theorem.

\bth\label{thm:Liouv} In $\mathrm{U}_q(\widehat{\gl}_{M|N})$, one has that
$$
B_q(L^{\pm}(zq^{-2}))=\z^{\pm}(zq^{2M-2N-2})B_q(L^{\pm}(z)).$$
\eth


Like quantum determinants \cite{KL}, we also find that the quantum Berezinian satisfies various minor identities such as
Jacobi's theorem (Theorem. \ref{Jac rat}), Schur complement theorem (Theorem. \ref{Sch com}) and Sylvester's theorem (Theorem. \ref{syl}) etc. However, there are subtitles in the super case, since these identities
 cannot be proved using similar methods as in quantum determinants. Instead we show these identities using the techniques of quasideterminants,
 so we also give decompositions of the quantum Berezinian into mutually commutative quasideterminants in Theorem \ref{thm:decomp}.

Finally we derive analogue of the MacMahon Master theorem \cite{MRa} for the quantum affine superalgebra $\mathrm{U}_q(\widehat{\gl}_{M|N})$.
\bth
We have the identity
\beq
\sum_{r=0}^k(-1)^r\str\mathcal{S}^q_r\mathcal{A}^q_{\{r+1,\ldots,k\}}L^{\pm}_1(z)L^{\pm}_2(zq^2)\cdots L^{\pm}_k(zq^{2k-2})=0,
\eeq
\beq
\sum_{r=0}^k(-1)^r\str\mathcal{A}^q_r\mathcal{S}^q_{\{r+1,\ldots,k\}}L^{\pm}_1(z)L^{\pm}_2(zq^2)\cdots L^{\pm}_k(zq^{2k-2})=0,
\eeq
where $\mathcal{A}^q_{\{r+1,\ldots,k\}}$ and $\mathcal{S}^q_{\{r+1,\ldots,k\}}$ denote the $q$-antisymmetrizer and $q$-symmetrizer over the copies of $\End (\BC^{M|N})$ labeled by $r+1,\cdots,k$.
\eth

We remark that many results in quantum affine superalgebras can be compared with those for the Yangians \cite{Mo1, Mo2, MRa} and
 quantum affine algebras \cite{JLM, L}.

\section{Quantum affine superalgebra $\mathrm{U}_q(\widehat{\gl}_{M|N})$}

In this section, we present the RLL realization of the quantum affine superalgebra $\mathrm{U}_q(\widehat{\gl}_{M|N})$ with $M\neq N$.
For simplicity we use the following standard parity: for $1\leq i\leq M+N$, 
\[
\bar{i}=\left\{\begin{array}{cc}
0&\text{ if }i\leq M\\
1&\text{ if }i> M,
\end{array}\right.
\]
Let $q_i=q^{1-2\bar i}$ for $i\in \{1,\ldots,M+N\}$, where $q$ is not a root of unity. Let $E_{ij}$ be the unit matrix in $\End \mathbb{C}^{M|N} $,
the $R$-matrix $R(z,w) \in (\End  \mathbb{C}^{M|N} \otimes \End  \mathbb{C}^{M|N} )[z,w]$ is defined as
\begin{eqnarray}
\begin{array}{rcl}
R(z,w) &=&  \sum\limits_{i\in I}(zq_i - wq_i^{-1}) E_{ii} \otimes E_{ii}  + (z-w) \sum\limits_{i \neq j} E_{ii} \otimes E_{jj} \\
&\ & + z\sum\limits_{i<j}(q_j-q_j^{-1})  E_{ij} \otimes E_{ji}+w \sum\limits_{i>j} (q_j-q_j^{-1}) E_{ij} \otimes E_{ji} .
\end{array}
\end{eqnarray}
It satisfies the following {\it quantum Yang-Baxter equation}
\begin{equation}\label{QYBE}
R_{12}(z_1,z_2)R_{13}(z_1,z_3)R_{23}(z_2,z_3) = R_{23}(z_2,z_3)R_{13}(z_1,z_3)R_{12}(z_1,z_2).
\end{equation}
Introduce the graded permutation operator $P$ on the tensor product $\mathbb C^{M|N}\otimes \mathbb C^{M|N}$ such that
$$P=\sum_{i,j=1}^{M+N}(-1)^{\bar{j}}\;E_{ij}\otimes E_{ji}.$$ Let $R = R(1,0),R' = P R^{-1}P$. Then $R(z,w) = z R - w R'$ and $R-R' =(q - q^{-1}) P$.
We also need the  $R$ matrix:
\beq
\ol R(z/w)=\frac{R(z,w)}{zq-wq^{-1}}.
\eeq
Denote $\ol R_{21}(z)=P_{12}\ol R_{12}(z)P_{12}$, then
\beql{Rin}
\ol R_{12}({z\over w})\ol R_{21}({w\over z})=1,\quad
\ol R_{21}({w\over z}) = \ol R_{q^{-1}}({z\over w}).
\eeq

For  a matrix $A=\sum _{i,j=1}^{M+N}a_{ij}E_{ij}$,   the supertranspose is defined by
$$A^{st}=\sum_{i,j}(-1)^{\bar{i}(\bar{i}+\bar{j})}a_{ji}E_{ij},$$
and the supertrace $\str$ by
\[
\str(A)=\sum_{i=1}^{M+N}(-1)^{\bar{i}}a_{ii}
\]
For any  $a\in\{1,2,\ldots,k\}$ we will denote by $st_a$ the corresponding partial transposition on the
algebra $(\End \mathbb C^{M|N})^{\otimes k}$ which acts as $st$ on the $a$-th copy of $\End \mathbb C^{M|N}$ and as the identity map on all
the other tensor factors.

We introduce the following normalized $R$-matrix
\beql{rf}
R(x)=f(x)\tss \overline R(x),
\eeq
where
\ben
f(x)=1+\sum_{k=1}^{\infty}f_kx^k,\qquad f_k=f_k(q),
\een
is a formal power series in $x$ 
uniquely determined by the relation
\ben
f(xq^{2N-2M})=f(x)\ts\frac{(1-xq^{-2})\tss(1-xq^{2N-2M+2})}{(1-x)\tss(1-xq^{2N-2M})}.
\een

The $R$-matrix \eqref{rf} satisfies the {\em crossing symmetry relations \cite{GZ}}:
\begin{align}\label{e:cross1}
\big(R_{12}(x)^{-1}\big)^{st_2} D_2 R_{12}(xq^{2N-2M})^{st_2}&=D_2,\\ \label{e:cross2}
R_{12}(xq^{2N-2M})^{st_1}\tss D_1\big(R_{12}(x)^{-1}\big)^{st_1}&=D_1,
\end{align}
where $D$ denotes the diagonal $n\times n$ matrix
\beql{d}
D=\diag\big[q^{2},q^{4},\ldots,q^{2M},q^{2M},\ldots,q^{2M-2N+2}\big].
\eeq

Following \cite{RS} we introduce the quantum affine superalgebra in the R-matrix format.
\defin
The quantum affine superalgebra $\mathrm U_q(\wh\gl_{M|N})$ is the superalgebra over $\CC(q)$,
generated by ${l^{\pm}_{ij}}^{(r)}$ and the invertible central element $q^c$, where $1\leq i, j\leq M+N$ and $r$ runs over nonnegative integers. Let $L^{\pm}(u)=(l_{ij}^{\pm}(u))$ be the matrix
\begin{equation}
\begin{split}
L^{\pm}(u)= \sum_{i,j=1}^{M+N}{  }l^{\pm}_{ij}(u)\otimes E_{ij},
\end{split}
\end{equation}
where $l^{\pm}_{ij}(u)$  are formal series in $u^{\mp 1}$ respectively:
\begin{equation}
\begin{split}
l^{\pm}_{ij}(u)= \sum_{r=0}^{\infty} {l^{\pm}_{ij}}^{(r)}u^{\mp r}.
\end{split}
\end{equation}
The defining relations are
\begin{align} \label{RLL-0}
{l^+_{ji}}^{(0)}={l^-_{ij}}^{(0)}=0,\ 1\leq i<j\leq M+N, \\ \label{RLL-1}
l{^-_{ii}}^{(0)} {l^+_{ii}}^{(0)}={l^+_{ii}}^{(0)}{l^-_{ii}}^{(0)}=1,\ 1\leq i\leq M+N, \\ \label{RLL}
R(z/w)L^{\pm}_{1} (z) L^{\pm}_{2} (w)= L^{\pm}_{2} (w) L^{\pm}_{1} (z)R(z/w),\\ \label{RLL cros}
R(z q^{c} /w) L^{+}_{1} (z) L^{-}_{2} (w)= L^{-}_{2} (w) L^{+}_{1} (z)R(z q^{-c}  /w).
\end{align}
The $q$-Yangian $\mathrm Y_q^+(\gl_{M|N})$ is the algebra generated by the $({l^+_{ij}}^{(0)})^{-1}$,
$ {l^+_{ij}}^{(r)}$ for $1\leq i,j \leq M+N  $ and $r \in \BZ_{\geq 0}$ with relations \eqref{RLL-0} and \eqref{RLL} \cite{Ch}.
The $q$-Yangian $\mathrm Y_q^-(\gl_{M|N})$ is defined similarly.

Note that when $N=0$ the quantum affine superalgebra $\mathrm U_q(\wh\gl_{M|0})$ descends to the quantum affine algebra
$\mathrm U_q(\wh\gl_{M})$ studied by Ding and Frenkel \cite{DF}.


\ble\label{RLL gen}
The defining relations \eqref{RLL} and \eqref{RLL cros} can be written in  terms of $l^{\pm}_{ij}(u)$ as
\beql{rll2}
\begin{split}
&\delta_{a=c}(zq_a-wq_a^{-1})l_{ab}^{\pm}(z)l_{cd}^{\pm}(w)(-1)^{(\bar{a}+\bar{b})(\bar{c}+\bar{d})}+\delta_{a\neq c}(z-w)l_{ab}^{\pm}(z)l_{cd}^{\pm}(w)(-1)^{(\bar{a}+\bar{b})(\bar{c}+\bar{d})}\\
&+\delta_{a>c}w(q-q^{-1})l_{cb}^{\pm}(z)l_{ad}^{\pm}(w)(-1)^{\bar{c}\bar{b}+\bar{c}\bar{d}+\bar{b}\bar{d}}
+\delta_{a<c}z(q-q^{-1})l_{cb}^{\pm}(z)l_{ad}^{\pm}(w)(-1)^{\bar{c}\bar{b}+\bar{c}\bar{d}+\bar{b}\bar{d}}\\
=&\delta_{b=d}(zq_b-wq_b^{-1})l_{cd}^{\pm}(w)l_{ab}^{\pm}(z)+\delta_{b\neq d}(z-w)l_{cd}^{\pm}(w)l_{ab}^{\pm}(z)\\
&+\delta_{d>b}w(q-q^{-1})l_{cb}^{\pm}(w)l_{ad}^{\pm}(z)(-1)^{\bar{c}\bar{b}+\bar{c}\bar{d}+\bar{b}\bar{d}}
+\delta_{d<b}z(q-q^{-1})l_{cb}^{\pm}(w)l_{ad}^{\pm}(z)(-1)^{\bar{c}\bar{b}+\bar{c}\bar{d}+\bar{b}\bar{d}}
\end{split}
\eeq
and
\beql{rll cros2}
\begin{split}
&\delta_{a=c}(zq_a^{c+1}-wq_a^{-1})l_{ab}^+(z)l_{cd}^-(w)(-1)^{(\bar{a}+\bar{b})(\bar{c}+\bar{d})}+\delta_{a\neq c}(zq^c-w)l_{ab}^+(z)l_{cd}^-(w)(-1)^{(\bar{a}+\bar{b})(\bar{c}+\bar{d})}\\
&+\delta_{a>c}w(q-q^{-1})l_{cb}^+(z)l_{ad}^-(w)(-1)^{\bar{c}\bar{b}+\bar{c}\bar{d}+\bar{b}\bar{d}}
+\delta_{a<c}zq^c(q-q^{-1})l_{cb}^+(z)l_{ad}^-(w)(-1)^{\bar{c}\bar{b}+\bar{c}\bar{d}+\bar{b}\bar{d}}\\
=&\delta_{b=d}(zq^{-c}q_b-wq_b^{-1})l_{cd}^-(w)l_{ab}^+(z)+\delta_{b\neq d}(zq^{-c}-w)l_{cd}^-(w)l_{ab}^+(z)\\
&+\delta_{d>b}w(q-q^{-1})l_{cb}^-(w)l_{ad}^+(z)(-1)^{\bar{c}\bar{b}+\bar{c}\bar{d}+\bar{b}\bar{d}}
+\delta_{d<b}zq^{-c}(q-q^{-1})l_{cb}^-(w)l_{ad}^+(z)(-1)^{\bar{c}\bar{b}+\bar{c}\bar{d}+\bar{b}\bar{d}}.
\end{split}
\eeq
\ele
\bpf These follow by applying both sides of \eqref{RLL} and \eqref{RLL cros} to $e_b\ot e_d$ and comparing the coefficient of $e_a\ot e_c$.
\epf

\ble\label{om}
The following mapping $\om_{M|N}$ defines an anti-automorphism of $\U_q(\wh\gl_{M|N})$
\beql{anti}
\om_{M|N}(L^{\pm}(z))= L^{\pm}(z)^{-1}.
\eeq
\ele
\bpf
The equation \eqref{RLL} is equivalent to
\beql{rtt opp}
R(z/w)L_2^{-1}(w)  L_1^{-1}(z) = L_1^{-1}(z) L_2^{-1}(w) R(z/w).
\eeq
Note that $\om_{M|N}$ is an involution. This completes the proof.
\epf

By relation \eqref{Rin},
multiplying both sides of \eqref{rtt opp} by the  inverse of $R(z/w)$ from the left and right we get that
\beql{rtt inv}
R_{q^{-1}}(z/w) L_1^{-1}(z)L_2^{-1}(w)  = L_2^{-1}(w) L_1^{-1}(z) R_{q^{-1}}(z/w).
\eeq

Denote the generator matrix of $\mathrm U_q(\wh\gl_{N|M})$ by $\bar L^{\pm}(z)$. We have the following isomorphism.
\ble\label{reve mn}
The mapping $\rho_{M|N}:\mathrm Y^{\pm}_q(\gl_{M|N})\rightarrow \mathrm Y^{\mp}_q(\gl_{N|M})$ defined by
\beq
\rho_{M|N}(l^{\pm}_{ij}(z))=\bar l^{\mp}_{M+N+1-i,M+N+1-j}(1/z)
\eeq
is an algebra isomorphism.
\ele
\bpf Due to the definition of the $q$-Yangian in \cite{Ch}, we only need to verify \eqref{RLL-1} and \eqref{RLL}. But the former one is clear.
From the RLL relation in terms of generators given in Lemma \ref{RLL gen}, we see that $\rho_{M|N}$ obeys the relation.
It is obviously that $\rho_{M|N}$ is bijective.
\epf

We will simply denote  $\om_{M|N}$ and $\rho_{M|N}$ by $\om $ and $\rho $ respectively if there is no confusion in the context.

There is a natural Hopf algebra structure on $\U_q(\wh\gl_{M|N})$ with the coassociative comultiplication homomorphism $\De:\U_q(\wh\gl_{M|N})\to\U_q(\wh\gl_{M|N})\ot\U_q(\wh\gl_{M|N})$ given by
\begin{equation}
\label{3.7}
l^{\pm}_{ij}(z)\mapsto\sum_{k=1}^{M+N}
l^{\pm}_{ik}(zq^{\pm1\ot\frac{c}{2}})\ot l^{\pm}_{kj}(zq^{\mp\frac{c}{2}\ot1})\,
{(-1)}^{\ts(\ts\bi\ts+\ts\bk\ts)(\ts\bj\ts+\ts\bk\ts)},
\end{equation}
and the counit homomorphism $\epsilon:\U_q(\wh\gl_{M|N})\to\CC$ is defined
by the assignment $L^{\pm}_{ij}(u)\mapsto\delta_{ij}\,$.
The antipodal mapping $S:\U_q(\wh\gl_{M|N})\to\U_q(\wh\gl_{M|N})$ is the antiautomorphism $\om_{M|N}$ defined by \eqref{anti}

\section{The center of $\mathrm{U}_q(\widehat{\gl}_{M|N})$ }
In the super Yangian $\mathrm Y(\gl_{M|N})$, Nazarov constructed two families of central elements \cite{Na}.
The first set was the coefficients of the Berezinian, and the second was related
to the Berezinian by the Liouville-type formulas.
In this section, we will construct two families of central elements
for $\mathrm{U}_q(\widehat{\gl}_{M|N})$ which are quantum analog of the central elements given by Nazarov.

The {\em Hecke algebra} $\Hc_m$ is generated by the elements
$T_1,\dots,T_{m-1}$ subject to the relations
\ben
\bal
&(T_i-q)(T_i+q^{-1})=0,\\
&T_{i}T_{i+1}T_{i}=T_{i+1}T_{i}T_{i+1},\\
&T_iT_j=T_jT_i \quad\text{for\ \  $|i-j|>1$}.
\eal
\een

For $i=1, \dots, m-1$
we let $\si_{i}=(i, i+1)$ be the adjacent transposition in the symmetric group $\Sym_{m}$.
Choose a reduced decomposition
$\si=\si_{i_{1}} \dots \si_{i_{l}}$ of any element
$\si \in \Sym_{m}$ and set $T_{\si}=T_{i_{1}} \dots T_{i_{l}}$. This element does not depend on
the choice of reduced decomposition of $\si$.

Setting $\check{R}=PR$, we get
\ben
\check{R}_{k}\check{R}_{k+1}\check{R}_{k}=\check{R}_{k+1}\check{R}_{k}\check{R}_{k+1}\Fand
(\check{R}_{k}-q)(\check{R}_{k}+q^{-1})=0,
\een
where $\check{R}_{k}=P_{k,k+1}R_{k,k+1}$.
Then we obtain a representation of the Hecke algebra $\Hc_{m}$ on the tensor product space $(\CC^{M|N})^{\ot m}$
defined by \cite{Jim}
\beql{haact}
T_{k}\mapsto \check{R}_{k},\qquad k=1,\dots,m-1.
\eeq

The $q$-symmetrizer $\mathcal S^q_m$ and $q$-antisymmetrizer $\mathcal A^q_m$ of the Hecke algebra are defined respectively as
\begin{align}
\mathcal S^q_{m} &=\frac{q^{-m(m-1)/2 }}{[m]_q!} \sum_{\sigma\in \mathfrak{S}_{m}}   q^{l(\si)}T_{\sigma},\\
\mathcal A^q_{m} &=\frac{q^{m(m-1)/2 }}{[m]_q!} \sum_{\sigma\in \mathfrak{S}_{m}}   (-q)^{-l(\si)}T_{\sigma},
\end{align}
wHere $[m]_q!=[1]_q\cdots [m]_q $ and $[n]_q=\frac{q^n-q^{-n}}{q-\qin}$.

Let $\hat{R}(x)=PR(x,1)$. Then the quantum Yang-Baxter equation \eqref{QYBE} implies that
\beql{QYBE2}
\hat{R}_{12}(z/w)\hat{R}_{23}(z/v)\hat{R}_{12}(w/v)=\hat{R}_{23}(z/w)\hat{R}_{12}(z/v)\hat{R}_{23}(w/v).
\eeq
And the $q$-(anti-)symmetrizer can be written recursively in the following form \cite{Jim}
\beql{idem}
\begin{split}
	&\mathcal S^q_2=\frac{q^{-1}}{q^2-q^{-2}}\hat{R}(q^2),\quad \mathcal S^q_{m+1}=\frac{q^{-m}}{q^{m+1}-q^{-m-1}} \mathcal S^q_m\hat{R}_{m,m+1}(q^{2m}) \mathcal S^q_m,\\
	&\mathcal A^q_2=\frac{q}{q^2-q^{-2}}\hat{R}(q^{-2}),\quad \mathcal A^q_{m+1}=\frac{q^m}{q^{m+1}-q^{-m-1}}\mathcal{A}^q_m\hat{R}_{m,m+1}(q^{-2m})\mathcal{A}^q_m.\\
\end{split}
\eeq

Let $L^{\pm}(z)^*=(L^{\pm}(z)^{-1})^{st}$. It follows from \eqref{RLL} that
\beql{RLLst}
\hat{R}_{\qin}(w/z)L^{\pm}_1(z)^*L^{\pm}_2(w)^*=L^{\pm}_1(w)^*L^{\pm}_2(z)^*\hat{R}_{\qin}(w/z).
\eeq

Consider the  tensor product $\mathrm{U}_{q}(\wh{\gl}_{M|N}) \otimes(\End\mathbb{C}^{M|N})^{\otimes m}$,
the $RLL$ relation and the Yang-Baxter equation imply that
\begin{align}\notag
&\mathcal{A}^q_mL^{\pm}_{1} (z)L^{\pm}_{2} (zq^2)\cdots L^{\pm}_{m}(zq^{2m-2})\\
& \qquad =L^{\pm}_{1}(zq^{2m-2})L^{\pm}_{2} (zq^{2m-4})\cdots  L^{\pm}_{m} (z)\mathcal{A}^q_m.\\ \notag
&\mathcal{S}^{\qin}_mL^{\pm}_{1}(z)^*L^{\pm}_{2}(zq^{-2})^* \cdots L^{\pm}_{m}(zq^{2-2m})^*\\
& \qquad =L^{\pm}_{1}(zq^{2-2m})^*L^{\pm}_{2} (zq^{4-2m})^*\cdots  L^{\pm}_{m} (z)^*\mathcal{S}^{\qin}_m.
\end{align}

Denote by $\mathcal{E}$ (resp. $\mathcal{O}$) the projector onto the even (resp. odd) subspace of $\BC^{M|N}$. Let $\mathcal{E}_a$ and $\mathcal{O}_a$ denote the projectors acting on the $a$-th factor of the tensor product.
The the {\it quantum Berezinian} of  $L^{\pm}(z) $ is defined as
\beql{def:ber}
\begin{split}
B_q(L^{\pm}(z) )=&
\str_{1,\ldots,M+N}\mathcal{E}_1\cdots\mathcal{E}_M\mathcal{O}_{M+1}\cdots\mathcal{O}_{M+N}\mathcal{A}^q_M
L_1^{\pm}(z)\cdots L_M^{\pm}(zq^{2M-2})\\
&\times (\mathcal{S}^{\qin}_N)_{M+1,\ldots,M+N}L_{M+1}^{\pm}(zq^{2M-2})^*\cdots L_{M+N}^{\pm}(zq^{2M-2N})^*,
\end{split}
\eeq
where $(\mathcal{S}^{\qin}_N)_{M+1,\ldots,M+N}$ is the $N$-th $\qin$-symmetrizer acting on the $M+1,\ldots,M+N$-th factors.

It can be verified that
\beql{ber}
\bal
B_q
&(L^{\pm}(z) )
=\sum_{\sigma\in S_M}(-q)^{-l(\sigma)}
L^{\pm}(z)_{ {\sigma(1)}, 1}\cdots L^{\pm}(zq^{2M-2})_{ {\sigma (M)},M}\\
&\times \sum_{\tau\in S_N} (-q)^{-l(\tau)}  \big({L^{\pm}(zq^{2M-2}) }^{-1}\big)_{M+1,M+\tau(1)}
\cdots \big({L^{\pm}(zq^{2M-2N})}^{-1}\big)_{M+N,{M+\tau(N)}}.\\
\eal
\eeq

In the following we will give a family of central elements in $\mathrm{U}_q(\widehat{\gl}_{M|N})$ and show that they are related to the quantum Berezinian.
\bpr\label{prop:equiz}
There exist a series $\z^+(z)$ in $z^{-1}$ and a series $\z^-(z)$ in $z$ with coefficients
in the quantum affine superalgebra $\mathrm{U}_q(\widehat{\gl}_{M|N})$ such that
\begin{align}\label{lltr}
	L^{\pm}(zq^{2N-2M})^{st}DL^{\pm}(z)^*&=\z^{\pm}(z)D \\
	\intertext{and}
	L^{\pm}(z)^*D^{-1} L^{\pm}(zq^{2N-2M})^{st}&=\z^{\pm}(z)D^{-1}.
	\label{secllrt}
\end{align}
Moreover, the coefficients of the series $\z^{\pm}(u)$ belong to the center of the algebra
$\U_q(\wh{\gl}_{M|N})$.
\epr

\bpf The proof is a super-analogue of Proposition 2.1 in \cite{JLM}. Multiply both sides of the 
relation of \eqref{RLL} by $L^{\pm}_2(w)^{-1}$ from the left and the right
and apply the transposition $st_2$ to get
\ben
R(z/w)^{st_2}L^{\pm}_2(w)^*L^{\pm}_1(z)
=L^{\pm}_1(z)L^{\pm}_2(w)^*R(z/w)^{st_2}
\een
and hence
\beql{rtrain}
\big(R(z/w)^{st_2}\big)^{-1}L^{\pm}_1(z)L^{\pm}_2(w)^*=
L^{\pm}_2(w)^*L^{\pm}_1(z)\big(R(z/w)^{st_2}\big)^{-1}.
\eeq
Using the crossing symmetry relation \eqref{e:cross1} to replace the $R$-matrix by
\ben
\big(R(z/w)^{st_2}\big)^{-1}=D_2^{-1}\big(R(q^{2M-2N} z/w)^{-1}\big)^{st_2}D_2,
\een
we get that
\begin{align*}
&\big(R(q^{2M-2N} z/w )^{-1}\big)^{st_2}D_2L^{\pm}_1(z)L^{\pm}_2(w)^*D_2^{-1}\\
&\qquad\qquad =D_2L^{\pm}_2(w)^*L^{\pm}_1(z)D_2^{-1}\big(R(q^{2M-2N} z/w )^{-1}\big)^{st_2}.
\end{align*}
Observe
that the $R$-matrix $R- R'$ equals $(q-\qin)P$, where $P$
is the permutation operator. Therefore,
\ben
\big((R- R')^{-1}\big)^{st_2}=\frac{1}{q-\qin}\ts Q\qquad\text{with}\quad
Q=\sum_{i,j=1}^{M+N} e_{ij}\ot e_{ij}(-1)^{\bar{i}\bar{j}+\bar i+\bar j}.
\een
Hence, by taking $z=wq^{2N-2M}$ we get
\beq\label{QLL}
Q D_2L^{\pm}_1(wq^{2N-2M})L^{\pm}_2(w)^*D_2^{-1}
=D_2L^{\pm}_2(w)^*L^{\pm}_1(wq^{2N-2M})D_2^{-1}Q.
\eeq
Since $Q$ is an operator in $\End\CC^{M|N}\ot\End\CC^{M|N}$ with a one-dimensional image,
both sides must be equal to $Q\tss \z^{\ts\pm}(w)$ for some series $\z^{\ts\pm}(w)$
with coefficients in the quantum affine superalgebra.
Using the relations $QX_1=QX_{2}^{st}$ and $X_1 Q=X_{2}^{st} Q$ which hold for an arbitrary matrix $X$, we can
write the equation of $\z^{\ts\pm}(w)$ as
\begin{align}\label{qld}
	QL^{\pm}_2(wq^{2N-2M})^{st}D_2L^{\pm}_2(w)^*&=Q\tss D_2\tss \z^{\ts\pm}(w)\\
	\intertext{and}
	L^{\pm}_2(w)^* D_2^{-1}L^{\pm}_2(wq^{2N-2M})^{st}Q&=D_2^{-1}\tss Q\tss \z^{\ts\pm}(w).
	\label{opqld}
\end{align}
By taking supertrace over the first copy of $\End\CC^{M|N}$ on both sides of \eqref{qld}
and \eqref{opqld} we arrive at
\eqref{lltr} and \eqref{secllrt}, respectively.

We will now use \eqref{secllrt}
to show that the series $\z^-(w)$ commutes with $L^+(z)$.
We have
\beql{loz}
L^+_1(z)\z^-(w)=L^+_1(z)D_2L^{-}_2(w)^* D_2^{-1}L^{-}_2(wq^{2N-2M})^{st}.
\eeq
Transform the relation \eqref{RLL cros} 
in the same way as we did for the relation \eqref{RLL}
in the beginning of the proof to get the following counterpart of \eqref{rtrain}:
\ben
\big(R(zq^c/w)^{st_2}\big)^{-1}L^{+}_1(z)L^{-}_2(w)^*=
L^{-}_2(w)^*L^{+}_1(z)\big(R(zq^{-c}/w)^{st_2}\big)^{-1}.
\een
Hence, the right hand side of \eqref{loz} equals
\beql{drt}
D_2 R(zq^{c}/w)^{st_2}L^{-}_2(w)^* L^+_1(z)
\big(R(zq^{-c}/w)^{st_2}\big)^{-1}D_2^{-1}L^{-}_2(wq^{2N-2M})^{st}.
\eeq
Applying again the crossing symmetry \eqref{e:cross1}, we can write
\ben
\big(R(zq^{-c}/w)^{st_2}\big)^{-1}D_2^{-1}=D_2^{-1}\big(R(zq^{-c}/wq^{2N-2M})^{-1}\big)^{st_2}.
\een
Continue transforming \eqref{drt}
by using the following consequence of \eqref{RLL cros}:
\beq
\begin{split}
&L^+_1(z)\big(R(zq^{-c}/wq^{2N-2M})^{-1}\big)^{st_2}L^{-}_2(wq^{2N-2M})^{st}\\
=&L^{-}_2(wq^{2N-2M})^{st}\big(R(zq^{c}/wq^{2N-2M})^{-1}\big)^{st_2}L^+_1(z),
\end{split}
\eeq
so that \eqref{drt} becomes
\ben
D_2 R(zq^{c}/w)^{st_2}L^{-}_2(w)^* D_2^{-1}L^{-}_2(wq^{2N-2M})^{st}
\big(R(zq^{c}/wq^{2N-2M})^{-1}\big)^{st_2} L^+_1(z).
\een
By \eqref{secllrt} this simplifies to
\ben
D_2 R(zq^{c}/w)^{st_2}D_2^{-1}
\big(R(zq^{c}/wq^{2N-2M})^{-1}\big)^{st_2}\z^-(w) L^+_1(z)
\een
which equals $\z^-(w)L^+_1(z)$ by \eqref{e:cross1}.
This proves that $L^+_1(z)\z^-(w)=\z^-(w)L^+_1(z)$. The commutation relation $L^-_1(z)\z^-(w)=\z^-(w)L^-_1(z)$
and the centrality of $\z^+(w)$ are verified in the same way.
\epf
Apply the relation \eqref{QLL} to $e_i\otimes e_j$ for some $1\leq j\leq M$ and compare the coefficients of $e_i\otimes e_j$ to get that
\[
\begin{split}
\de_{ij}\z^{\pm}(w)=&\sum_{k=1}^Ml_{ki}^{\pm}(wq^{2N-2M})l^{\pm'}_{jk}(w)q^{2k-2j}-\sum_{k=M+1}^{M+N}l_{ki}^{\pm}(wq^{2N-2M})l^{\pm'}_{jk}(w)q^{4M+2-2k-2j}.\\
\end{split}
\]
While for $M+1\leq j\leq M+N$, there is
\[
\begin{split}
\de_{ij}\z^{\pm}(w)=&\sum_{k=1}^Ml_{ki}^{\pm}(wq^{2N-2M})l^{\pm'}_{jk}(w)q^{2k+2j-4M-2}-\sum_{k=M+1}^{M+N}l_{ki}^{\pm}(wq^{2N-2M})l^{\pm'}_{jk}(w)q^{-2k+2j}.
\end{split}
\]
Here, we use $l^{\pm'}_{ij}(w)$ to denote the entries in $L^{\pm}(w)^{-1}$.
By applying the square of antipode $S^2$ to $L^{\pm}(w)L^{\pm}(w)^{-1}=1$, we can see that
\beql{square anti}
\begin{split}
\z^{\pm}(w)S^2(l^{\pm}_{ij}(w))=\begin{cases} &q^{2i-2j}l^{\pm}_{ij}(w) \text{ when } 1\leq i,j\leq M;\\
&q^{4M+2-2i-2j}l^{\pm}_{ij}(w)\text{ when } 1\leq j\leq M, M+1\leq i\leq M+N;\\
&q^{2i+2j-4M-2}l^{\pm}_{ij}(w)\text{ when } 1\leq i\leq M, M+1\leq j\leq M+N;\\
&q^{-2i+2j}l^{\pm}_{ij}(w)\text{ when } M+1\leq i,j\leq M+N.\\
\end{cases}
\end{split}
\eeq
\bco\label{cop on z}
$\Delta(\z^{\pm}(w))=\z^{\pm}(wq^{\pm1\ot\frac{c}{2}})\ot \z^{\pm}(wq^{\mp\frac{c}{2}\ot1})$.
\eco
\bpf
It is a direct result of the fact that the square of antipode is a coalgebra homomrphism, i.e. $\Delta S^2=(S^2\ot S^2)\Delta$. By applying both sides of the equality to $l^{\pm}_{ij}(w)$ and using \eqref{square anti}, we get the result.
\epf

\bco\label{cor:zu}
If $M\neq N$, we have the formulas
\begin{align}\label{zu}
	\z^{\pm}(z)&=\frac{1 }{q^{M+1}[M]_q-q^{2M-N+1}[N]_q}\ts \str\ts DL^{\pm}(zq^{2N-2M})L^{\pm}(z)^{-1}\\
	\intertext{and}
	\z^{\pm}(z)&=\frac{1}{q^{-(M+1)}[M]_q-q^{-2M+N-1}[N]_q}\ts \str\ts D^{-1}L^{\pm}(z)^{-1}L^{\pm}(zq^{2N-2M}).
	\label{opzu}
\end{align}
\eco

\bpf
Using the relation
\beq
\str(XY)=\str(X^{st}Y^{st}),
\eeq
the formulas follow by taking super trace on both sides of the respective matrix relations
\eqref{lltr} and \eqref{secllrt}.
\epf

The center elements $\z^{\pm}(z)$ are related to the quantum Berezinian by the following  Liouville formula.

\bth\label{thm:Liouv} In $\mathrm{U}_q(\widehat{\gl}_{M|N})$, one has that
$$
B_q(L^{\pm}(zq^{-2}))=\z^{\pm}(zq^{2M-2N-2})B_q(L^{\pm}(z)).$$
\eth
\bpf
%
We use the similar appraoch of the proof of Theorem 4.1 in \cite{Na2}.
It follows from \eqref{RLL} and \eqref{RLLst} that
\beql{RLLL}
\begin{split}
&\mathcal{R}L_1^{\pm}(zq^{-2})L_2^{\pm}(z)\cdots L_{M+1}^{\pm}(zq^{2M-2})\\
&\times L_{M+2}^{\pm}(zq^{2M-2})^*L_{M+3}^{\pm}(zq^{2M-4})^*\cdots L_{M+N+2}^{\pm}(zq^{2M-2N-2})^*\\
=&L_1^{\pm}(z)\cdots L_M^{\pm}(zq^{2M-2})L_{M+1}^{\pm}(zq^{-2})\\
&\times L_{M+2}^{\pm}(zq^{2M-2N-2})^*L_{M+3}^{\pm}(zq^{2M-2})^*\cdots L_{M+N+2}^{\pm}(zq^{2M-2N})^*\mathcal{R},
\end{split}
\eeq
where
\[
\begin{split}
\mathcal{R}=&\hat{R}_{M,M+1}(q^{-2M})\hat{R}_{M-1,M}(q^{2-2M})\cdots \hat{R}_{1,2}(q^{-2})\\
&\times (\hat{R}_{\qin})_{M+2,M+3}(q^{-2N})(\hat{R}_{\qin})_{M+3,M+4}(q^{2-2N})\cdots (\hat{R}_{\qin})_{M+N+1,M+N+2}(q^{-2}).
\end{split}
\]

Multiply both sides of \eqref{RLLL} by
\[
\begin{split}
&\mathcal{E}_1\cdots\mathcal{E}_{M}\mathcal{O}_{M+3}\cdots\mathcal{O}_{M+N+2}(\mathcal{A}_{M-1}^q)_{1,\ldots,M-1}\times(\mathcal{S}_{N-1}^{\qin})_{M+4,\ldots, M+N+2}Q_{M+1,M+2}D_{M+2}
\end{split}
\]
on the left and by $D_{M+2}^{-1}$ on the right. Then the left side of the equation equals to
\beql{ls}
\begin{split}
&\mathcal{E}_1\cdots\mathcal{E}_{M}\mathcal{O}_{M+3}\cdots\mathcal{O}_{M+N+2}(\mathcal{A}_{M-1}^q)_{1,\ldots,M-1} (\mathcal{S}_{N-1}^{\qin})_{M+4,\ldots, M+N+2}Q_{M+1,M+2}D_{M+2}\\
&\times\hat{R}_{M,M+1}(q^{-2M})\hat{R}_{M-1,M}(q^{2-2M})\cdots \hat{R}_{1,2}(q^{-2})\\
&\times (\hat{R}_{\qin})_{M+2,M+3}(q^{-2N})(\hat{R}_{\qin})_{M+3,M+4}(q^{2-2N})\cdots (\hat{R}_{\qin})_{M+N+1,M+N+2}(q^{-2})\\
&\times L_1^{\pm}(zq^{-2})L_2^{\pm}(z)\cdots L_{M+1}^{\pm}(zq^{2M-2})\\
&\times L_{M+2}^{\pm}(zq^{2M-2})^*L_{M+3}^{\pm}(zq^{2M-4})^*\cdots L_{M+N+2}^{\pm}(zq^{2M-2N-2})^*D_{M+2}^{-1}.
\end{split}
\eeq
And the right side is
\beql{rs}
\begin{split}
&\mathcal{E}_1\cdots\mathcal{E}_{M}\mathcal{O}_{M+3}\cdots\mathcal{O}_{M+N+2}(\mathcal{A}_{M-1}^q)_{1,\ldots,M-1} (\mathcal{S}_{N-1}^{\qin})_{M+4,\ldots, M+N+2}Q_{M+1,M+2}D_{M+2}\\
&L_1^{\pm}(z)\cdots L_M^{\pm}(zq^{2M-2})L_{M+1}^{\pm}(zq^{-2})\\
&\times L_{M+2}^{\pm}(zq^{2M-2N-2})^*L_{M+3}^{\pm}(zq^{2M-2})^*\cdots L_{M+N+2}^{\pm}(zq^{2M-2N})^*\\
&\times\hat{R}_{M,M+1}(q^{-2M})\hat{R}_{M-1,M}(q^{2-2M})\cdots \hat{R}_{1,2}(q^{-2})\\
&\times (\hat{R}_{\qin})_{M+2,M+3}(q^{-2N})(\hat{R}_{\qin})_{M+3,M+4}(q^{2-2N})\cdots (\hat{R}_{\qin})_{M+N+1,M+N+2}(q^{-2})D_{M+2}^{-1}.
\end{split}
\eeq
By \eqref{idem} and $\mathcal{A}_M^q=\mathcal{A}_M^q\mathcal{A}_{M-1}^q$, we have that
\beql{anti-symm}
(\mathcal{A}_{M-1}^q)_{1,\ldots,M-1} \hat{R}_{M-1,M}(q^{2-2M})\cdots \hat{R}_{1,2}(q^{-2})=a_M(\mathcal{A}_M^q)_{1,\ldots,M},
\eeq
for some nonzero $a_M\in \BC(q)$.
Similarly,
\beql{symm}
\begin{split}
&(\mathcal{S}_{N-1}^{\qin})_{M+4,\ldots, M+N+2}(\hat{R}_{\qin})_{M+3,M+4}(q^{2-2N})\cdots (\hat{R}_{\qin})_{M+N+1,M+N+2}(q^{-2})\\
=&s_N(\mathcal{S}_N^{\qin})_{M+3,\ldots, M+N+2}
\end{split}
\eeq
for some nonzero $s_N\in \BC(q)$.
Thus up to a nonzero scalar in $\BC(q)$, \eqref{ls} equals to
\beql{ls2}
\begin{split}
&\mathcal{E}_1\cdots\mathcal{E}_{M}\mathcal{O}_{M+3}\cdots\mathcal{O}_{M+N+2}Q_{M+1,M+2}D_{M+2}\\
&\times \hat{R}_{M,M+1}(q^{-2M})(\hat{R}_{\qin})_{M+2,M+3}(q^{-2N})(\mathcal{A}_M^q)_{1,\ldots,M}(\mathcal{S}_{N}^{\qin})_{M+3,\ldots, M+N+2}\\
&\times L_1^{\pm}(zq^{-2})L_2^{\pm}(z)\cdots L_{M+1}^{\pm}(zq^{2M-2})\\
&\times L_{M+2}^{\pm}(zq^{2M-2})^*L_{M+3}^{\pm}(zq^{2M-4})^*\cdots L_{M+N+2}^{\pm}(zq^{2M-2N-2})^*D_{M+2}^{-1}.
\end{split}
\eeq

By the RLL relations and the definitions of
$\mathcal{A}_M^q$ and $\mathcal{S}_{N}^{\qin}$, we have that
\[
\begin{split}
&(\mathcal{A}_M^q)_{1,\ldots,M}L_1^{\pm}(zq^{-2})\cdots L_{M+1}^{\pm}(zq^{2M-2})\\
=&(\mathcal{A}_M^q)_{1,\ldots,M}L_1^{\pm}(zq^{-2})\cdots L_{M+1}^{\pm}(zq^{2M-2})(\mathcal{A}_M^q)_{1,\ldots,M},
\end{split}
\]
and
\[
\begin{split}
&(\mathcal{S}_{N}^{\qin})_{M+3,\ldots, M+N+2}L_{M+2}^{\pm}(zq^{2M-2})^*\cdots L_{M+N+2}^{\pm}(zq^{2M-2N-2})^*\\
=&(\mathcal{S}_{N}^{\qin})_{M+3,\ldots, M+N+2}L_{M+2}^{\pm}(zq^{2M-2})^*\cdots L_{M+N+2}^{\pm}(zq^{2M-2N-2})^*(\mathcal{S}_{N}^{\qin})_{M+3,\ldots, M+N+2}.
\end{split}
\]
Applying the operations $\str_{1},\ldots, \str_M$, $\str_{M+3}, \ldots,\str_{M+N+2}$ to \eqref{ls2} and notice that for homogenous elements $X,X'\in \End \BC^{M|N}$, $\str(XX')=\str(X'X)(-1)^{\deg X\deg X'}$, we get that
\beql{ls3}
\begin{split}
&\str_{1}\cdots \str_M\str_{M+3} \cdots \str_{M+N+2}\mathcal{E}_1\cdots\mathcal{E}_{M}\mathcal{O}_{M+3}\cdots\mathcal{O}_{M+N+2}Q_{M+1,M+2}D_{M+2}\\
&\times \hat{R}_{M,M+1}(q^{-2M})(\hat{R}_{\qin})_{M+2,M+3}(q^{-2N})(\mathcal{A}_M^q)_{1,\ldots,M}(\mathcal{S}_{N}^{\qin})_{M+3,\ldots, M+N+2}\\
&\times L_1^{\pm}(zq^{-2})L_2^{\pm}(z)\cdots L_{M+1}^{\pm}(zq^{2M-2})\\
&\times L_{M+2}^{\pm}(zq^{2M-2})^*L_{M+3}^{\pm}(zq^{2M-4})^*\cdots L_{M+N+2}^{\pm}(zq^{2M-2N-2})^*D_{M+2}^{-1}\\
=&\str_{1}\cdots \str_M\str_{M+3} \cdots \str_{M+N+2}\mathcal{E}_1\cdots\mathcal{E}_{M}\mathcal{O}_{M+3}\cdots\mathcal{O}_{M+N+2}Q_{M+1,M+2}D_{M+2}\\
&\times \hat{R}_{M,M+1}(q^{-2M})(\hat{R}_{\qin})_{M+2,M+3}(q^{-2N})(\mathcal{A}_M^q)_{1,\ldots,M}(\mathcal{S}_{N}^{\qin})_{M+3,\ldots, M+N+2}\\
&\times L_1^{\pm}(zq^{-2})L_2^{\pm}(z)\cdots L_{M+1}^{\pm}(zq^{2M-2})\\
&\times L_{M+2}^{\pm}(zq^{2M-2})^*\cdots L_{M+N+2}^{\pm}(zq^{2M-2N-2})^*D_{M+2}^{-1}(\mathcal{A}_M^q)_{1,\ldots,M}(\mathcal{S}_{N}^{\qin})_{M+3,\ldots, M+N+2}\\
=&\str_{1}\cdots \str_M\str_{M+3} \cdots \str_{M+N+2}\mathcal{E}_1\cdots\mathcal{E}_{M}\mathcal{O}_{M+3}\cdots\mathcal{O}_{M+N+2}Q_{M+1,M+2}D_{M+2}\\
&\times(\mathcal{A}_M^q)_{1,\ldots,M}(\mathcal{S}_{N}^{\qin})_{M+3,\ldots, M+N+2}\\
&\times \hat{R}_{M,M+1}(q^{-2M})(\hat{R}_{\qin})_{M+2,M+3}(q^{-2N})(\mathcal{A}_M^q)_{1,\ldots,M}(\mathcal{S}_{N}^{\qin})_{M+3,\ldots, M+N+2}\\
&\times L_1^{\pm}(zq^{-2})L_2^{\pm}(z)\cdots L_{M+1}^{\pm}(zq^{2M-2})\\
&\times L_{M+2}^{\pm}(zq^{2M-2})^*\cdots L_{M+N+2}^{\pm}(zq^{2M-2N-2})^*D_{M+2}^{-1}\\
=&\str_{1}\cdots \str_M\str_{M+3} \cdots \str_{M+N+2}\mathcal{E}_1\cdots\mathcal{E}_{M}\mathcal{O}_{M+3}\cdots\mathcal{O}_{M+N+2}\\
&\times c(\mathcal{A}_{M+1}^q)_{1,\ldots,M+1}(\mathcal{S}_{N+1}^{\qin})_{M+2,\ldots, M+N+2}\\
&\times L_1^{\pm}(zq^{-2})L_2^{\pm}(z)\cdots L_{M+1}^{\pm}(zq^{2M-2})\\
&\times L_{M+2}^{\pm}(zq^{2M-2})^*\cdots L_{M+N+2}^{\pm}(zq^{2M-2N-2})^*D_{M+2}^{-1}Q_{M+1,M+2}D_{M+2}\\
\end{split}
\eeq
for some nonzero $c\in \BC(q)$.
Notice that
\[
L_{M+1}^{\pm}(zq^{2M-2})L_{M+2}^{\pm}(zq^{2M-2})^*D_{M+2}^{-1}Q_{M+1,M+2}D_{M+2}=Q_{M+1,M+2}.
\]
Thus we can get that after applying $\str_{1}\cdots \str_M\str_{M+3} \cdots \str_{M+N+2}$ to \eqref{ls}, the result is $QB_q(L^{\pm}(zq^{-2}))$ up to a nonzero scalar.

For the same reason, up to a scalar in $\BC(q)$, \eqref{rs} equals to
\beql{rs2}
\begin{split}
&\mathcal{E}_1\cdots\mathcal{E}_{M}\mathcal{O}_{M+3}\cdots\mathcal{O}_{M+N+2}(\mathcal{A}_{M-1}^q)_{1,\ldots,M-1} (\mathcal{S}_{N-1}^{\qin})_{M+4,\ldots, M+N+2}Q_{M+1,M+2}D_{M+2}\\
&L_1^{\pm}(z)\cdots L_M^{\pm}(zq^{2M-2})L_{M+1}^{\pm}(zq^{-2})\\
&\times L_{M+2}^{\pm}(zq^{2M-2N-2})^*L_{M+3}^{\pm}(zq^{2M-2})^*\cdots L_{M+N+2}^{\pm}(zq^{2M-2N})^*\\
&\times\hat{R}_{M,M+1}(q^{-2M}) (\hat{R}_{\qin})_{M+2,M+3}(q^{-2N})(\mathcal{A}_{M}^q)_{1,\ldots,M} (\mathcal{S}_{N}^{\qin})_{M+3,\ldots, M+N+2}D_{M+2}^{-1}.
\end{split}
\eeq
Due to \eqref{QLL}, it equals to
\beql{rs3}
\begin{split}
&\mathcal{E}_1\cdots\mathcal{E}_{M}\mathcal{O}_{M+3}\cdots\mathcal{O}_{M+N+2}(\mathcal{A}_{M-1}^q)_{1,\ldots,M-1} (\mathcal{S}_{N-1}^{\qin})_{M+4,\ldots, M+N+2}Q_{M+1,M+2}\\
&L_1^{\pm}(z)\cdots L_M^{\pm}(zq^{2M-2}) L_{M+3}^{\pm}(zq^{2M-2})^*\cdots L_{M+N+2}^{\pm}(zq^{2M-2N})^*D_{M+2}\hat{R}_{M,M+1}(q^{-2M})\\
&\times  (\hat{R}_{\qin})_{M+2,M+3}(q^{-2N})(\mathcal{A}_{M}^q)_{1,\ldots,M} (\mathcal{S}_{N}^{\qin})_{M+3,\ldots, M+N+2}D_{M+2}^{-1}\z^{\pm}(zq^{2M-2N-2}).
\end{split}
\eeq
Applying $\str_{1}\cdots \str_M \str_{M+3}\cdots\str_{M+N+2}$ to it, we get that
\[
\begin{split}
&\str_{1}\cdots \str_M\str_{M+3} \cdots \str_{M+N+2}\mathcal{E}_1\cdots\mathcal{E}_{M}\mathcal{O}_{M+3}\cdots\mathcal{O}_{M+N+2}\\
&\times(\mathcal{A}_{M}^q)_{1,\ldots,M} (\mathcal{S}_{N}^{\qin})_{M+3,\ldots, M+N+2}Q_{M+1,M+2}\\
&L_1^{\pm}(z)\cdots L_M^{\pm}(zq^{2M-2}) L_{M+3}^{\pm}(zq^{2M-2})^*\cdots L_{M+N+2}^{\pm}(zq^{2M-2N})^*\\
&\times D_{M+2}\hat{R}_{M,M+1}(q^{-2M}) (\hat{R}_{\qin})_{M+2,M+3}(q^{-2N})D_{M+2}^{-1}\z^{\pm}(zq^{2M-2N-2}),
\end{split}
\]
which is $Q\z^{\pm}(zq^{2M-2N-2})B_q(L^{\pm}(z))$ up to a scalar. Since the constant term of $\z^{\pm}(z)$ is $1$, we get that
\[
B_q(L^{\pm}(zq^{-2}))=\z^{\pm}(zq^{2M-2N-2})B_q(L^{\pm}(z)).
\]
\epf

\bco
The coefficients of the quantum Berezinian belong to the center of  $\U_q(\wh{\gl}_{M|N})$ .
\eco
\bpf

The quantum Berezinian can be written in a formal power series
\beq
B_q(L^{\pm}(z))=b_0^{\pm}+b_1^{\pm}u^{\mp1}+b_2^{\pm}u^{\mp2}+\cdots.
\eeq
By  Theorem \ref{thm:Liouv}, it is sufficient to show that
$b_0^{\pm}$ are central elements.
It's easy to see that
\beq
b_0^{\pm}=l_{11}^{\pm(0)}l_{22}^{\pm(0)}\cdots l_{MM}^{\pm(0)}(l_{M+1,M+1}^{\pm(0)})^{-1}\cdots (l_{M+N,M+N}^{\pm(0)})^{-1}.
\eeq

Comparing coefficients of $w$ in the relation \eqref{rll2}, we have that

\beq
\bal
&l_{ii}^+(z)l_{kk}^{+(0)}=l_{kk}^{+(0)}l_{ii}^+(z),\qquad &1\leq i, k\leq M+N,\\
&l_{ij}^+(z)l_{kk}^{+(0)}=l_{kk}^{+(0)}l_{ij}^+(z) ,\qquad & k\neq i,j,\\
&l_{ij}^+(z)l_{ii}^{+(0)}l_{jj}^{+(0)}=l_{ii}^{+(0)}l_{jj}^{+(0)}l_{ij}^+(z),\qquad & \bar{i}=\bar{j},\\
&l_{ij}^+(z)l_{ii}^{+(0)}(l_{jj}^{+(0)})^{-1}=l_{ii}^{+(0)}(l_{jj}^{+(0)})^{-1}l_{ij}^+(z),\qquad &\bar{i}\neq\bar{j}.
\eal
\eeq
Thus, $b_0^{+}$ is central. Using the same arguments we can show that  $b_0^{-}$ is central.
%
%
%
\epf
Due to Corollary \ref{cop on z} and Theorem \ref{thm:Liouv}, we get that
\bco\label{cop on b}
$\Delta(B_q(L^{\pm}(z)))=B_q(L^{\pm}(z))\ot B_q(L^{\pm}(z)).$
\eco

\bre
When it comes to the $q$-Yangian $Y^{+}_q(\gl_{M|N})$, we can use the $R$-matrix $\bar{R}(x)$ to define the RLL relation. Then in the situation of $M=N$, it make sense and we readily get that all the above result holds.
\ere
%
%

\section{Minor identities for Berezinian}

Let $A$ be any square matrix of size $M+N$.
For any $i_1,\ldots,i_k,j_1,\ldots,j_k\in [1,M+N]$,
we denote by $A^{i_1,\ldots,i_k}_{j_1,\ldots,j_k}$ the matrix whose
$ab$-th entry is $A_{i_a j_b}$.
For any subsets $I=\{i_1<i_2<\cdots i_{k}\}$,
we denote by $A_I$ the submatrix of $A$ with rows and columns indexed by $I$.
In particular, we denote by $A^{(k)}$(resp. $A_{(k)}$) the principal submatrix of
$A$ of the form $(A_{ij})_{1\leq i,j\leq k}$ (resp. $(A_{ij})_{M+N-k+1\leq i,j\leq M+N}$ ).

For an invertible matrix $A$ with possible non-commutative entries, the quasideterminants $|A|_{i, j}$ is defined to be
$((A^{-1})_{j, i})^{-1}$. It also equals
\ben
\left|\begin{matrix}a_{11}&\dots&a_{1j}&\dots&a_{1N}\\
	&\dots&      &\dots&      \\
	a_{i1}&\dots&\boxed{a_{ij}}&\dots&a_{iN}\\
	&\dots&      &\dots&      \\
	a_{N1}&\dots&a_{Nj}&\dots&a_{NN}
\end{matrix}\right|=a_{ij}-r_i^j(A^{1,2,\ldots,\hat{i},\ldots,M+N}_{1,2,\ldots,\hat{j},\ldots,M+N})^{-1}c^i_j,
\een
where $r_{i}^j$ ($c_j^i$, respectively) is the row (column, respectively) matrix obtained from the $i$-th row ($j$-th column, respectively) of $A$ by deleting the element $a_{ij}$, and the symbol $\hat{}$ means the indice should be omitted (see \cite{GR}).

Let us recall the Jacobi ratio theorem from \cite{KL}
\bth[Gelfand and Retakh]
Let $A$ be the generic matrix of order $n$, let $B$ be its inverse and let $(\{i\},L,P)$ and $(\{j\},M,Q)$ be two partitions of $\{1,2,\ldots,n\}$ such that $|L|=|M|$ and $|P|=|Q|$. Then there holds:
\[
\left|B_{M\cup\{j\},L\cup\{i\}} \right|_{ji}=\left|A_{P\cup\{i\},Q\cup\{j\}} \right|_{ij}^{-1}.
\]
\eth

\bth\label{thm:decomp}
We have the following decomposition of $B_q(L^{\pm}(z) )$ in  $\mathrm{U}_q(\widehat{\gl}_{M|N})$,
\begin{multline}\label{decomp}
	B_q(L^{\pm}(z) )=\big|L^{\pm}(z)^{(1)}\big|_{11}\cdots
	\big|L^{\pm}(zq^{2M-2})^{(M)}\big|_{MM}\\[1em]
	\times\big|L^{\pm}(zq^{2M-2})^{(M+1)}\big|^{-1}_{M+1,M+1}\cdots
	\big|L^{\pm}(zq^{2M-2N})^{(M+N)}\big|^{-1}_{M+N,M+N}
\end{multline}
Moreover, the factors in the decomposition
are mutually commutative.
\eth

\bpf
The proof is similar to the case of super Yangian $Y(\gl_{m|n})$ given in \cite{G}.
We first note that in the quantum affine algebra $U_q(\wh{\gl}_M)$,
\beql{qdet}
\bal
{\det}_q(L^{\pm}(z)^{(M)})=&\sum_{\sigma\in S_M}(-q)^{-l(\sigma)}
L^{\pm}(z)_{ {\sigma(1)}, 1}\cdots L^{\pm}(zq^{2M-2})_{ {\sigma (M)},M}\\
=&\big|L^{\pm}(z)^{(1)}\big|_{11}\cdots
\big|L^{\pm}(zq^{2M-2})^{(M)}\big|_{MM}
\eal
\eeq
Thus, we get the first factors of \eqref{ber} and \eqref{decomp} coincide.

Let $\bar L^{\pm}(z)$ be the generator matrix of $\mathrm U_q(\wh\gl_{N|M})$,
then ${\bar L^{\pm}(z)}^{-1}$ satisfies the $\qin$-RLL relation and the quantum determinant
${\det}_{\qin}\left(( {\bar L}^{\pm}(q^{2N-2M}/z))^{-1})^{(N)}\right)$  can be expressed as
\beql{qdet2}
\begin{split}
	&\sum_{\tau\in S_N}(-q)^{-l(\tau)}
	({\bar L^{\pm}(q^{2-2M}/z)}^{-1})_{N,\tau(N)}\cdots ({\bar L^{\pm}(q^{2N-2M}/z)}^{-1})_{1,\tau (1)}\\
	=&\big|( {\bar L^{\pm}(q^{2N-2M}/z)}^{-1})^{(1)}\big|_{11}
	\cdots\big|({\bar L^{\pm}(q^{2-2M}/z)}^{-1})^{(N)}\big|_{N,N}.
\end{split}
\eeq
Apply  $\rho_{N|M}$ to Equation \eqref{qdet2} to get that
\beql{qdet3}
\begin{split}
	&\sum_{\tau\in S_k}(-q)^{-l(\tau)}(L^{\pm}(zq^{2M-2})^{-1})_{M+1,M+\tau(1)}\cdots (L^{\pm}(zq^{2M-2N})^{-1})_{M+N,M+\tau (N)}\\
	=&\big|(L^{\pm}(zq^{2M-2})^{-1})_{\{M+1,\ldots,M+N\}}\big|_{M+1,M+1}\cdots\big|(L^{\pm}(zq^{2M-2N})^{-1})_{\{M+N\}}\big|_{M+N,M+N}.
\end{split}
\eeq
By Jacobi's ratio theorem for quasideterminant \cite{KL},
\beq
\big|(L^{\pm}(zq^{2M-2k})^{-1})_{\{M+k,\ldots,M+N\}}\big|_{M+k,M+k}=\big|L^{\pm}(zq^{2M-2k})^{(M+k)}\big|^{-1}_{M+k,M+k}.
\eeq
Then equation \eqref{decomp} follows.

The factors in the decomposition
are mutually commutative since the coefficients of Berezinian belong to the center.
\epf

For any $k\times k$ matrix, we denote by $\pi (A)$ the matrix whose $(i,j)$-th entry is $A_{k+1-i,k+1-j}$.
The following theorem is an analog of Jacobi's ratio theorem for the quantum Berezinian.
\bth\label{Jac rat}
Let $I=(i_1<\cdots<i_k)$ be a subset of $[1,M+N]$  and
$I^{c}=\{i_{k+1}<\cdots<i_{M+N}\}$ be its complement in $[1,M+N]$.
If $I\subseteq [1,M]$ or $[1,M]\subseteq I $,
then
\beq
B_q(L^{\pm}(z))= B_q(L^{\pm}(z)_I) B_{\qin} \left(\pi((L^{\pm}(zq^{2M-2N})^{-1})_{I^{c}})\right).
\eeq
\eth

\bpf
We first prove the theorem in the case $I=(i_1<\cdots<i_k) \subseteq[1,M]$. In fact, due to Theorem 3.1 in \cite{KL}, we know that \eqref{qdet} is invariant under permutations of $\{1,2,\ldots,N\}$.
Thus, according to Theorem \ref{thm:decomp},
\beq
\begin{split}
B_q(L^{\pm}(z))
=&\big|L^{\pm}(z)^{i_1}_{i_1}\big|_{i_1i_1}
\big|L^{\pm}(zq^{2})^{i_1 i_2}_{i_1i_2}\big|_{i_2i_2}
\cdots|L^{\pm}(zq^{2M-2})^{i_1,\ldots,i_M}_{i_1,\ldots,i_M}\big|_{i_Mi_M}\\
&\times\big|(L^{\pm}(zq^{2M-2}))^{(M+1)}\big|^{-1}_{M+1M+1}\ldots\big|(L^{\pm}(zq^{2M-2N}))^{(M+N)}\big|^{-1}_{M+NM+N}\\
=&B_q(L^{\pm}(z)_I) \big|L^{\pm}(zq^{2k})^{i_1,\ldots,i_{k+1}}_{i_1,\ldots,i_{k+1}}\big|_{i_{k+1},i_{k+1}}
\cdots|L^{\pm}(zq^{2M-2})^{i_1,\ldots,i_M}_{i_1,\ldots,i_M}\big|_{i_Mi_M}\\
&\times\big|(L^{\pm}(zq^{2M-2}))^{(M+1)}\big|^{-1}_{M+1M+1}\ldots\big|(L^{\pm}(zq^{2M-2N}))^{(M+N)}\big|^{-1}_{M+NM+N}\\
\end{split}
\eeq
It follows from Jacobi's ratio theorem for quasideterminant that
\beq
\begin{split}
&\big|L^{\pm}(zq^{2k})^{i_1,\ldots,i_{k+1}}_{i_1,\ldots,i_{k+1}}\big|_{i_{k+1},i_{k+1}}
\cdots|L^{\pm}(zq^{2M-2})^{i_1,\ldots,i_M}_{i_1,\ldots,i_M}\big|_{i_Mi_M}\\
&\times\big|(L^{\pm}(zq^{2M-2}))^{(M+1)}\big|^{-1}_{M+1M+1}\ldots\big|(L^{\pm}(zq^{2M-2N}))^{(M+N)}\big|^{-1}_{M+NM+N}\\
=
&\big|(L^{\pm}(zq^{2k})^{-1})^{i_{k+1},\ldots,i_{M+N}}_{i_{k+1},\ldots,i_{M+N}}\big|_{i_{k+1}i_{k+1}}^{-1}\cdots
\big|(L^{\pm}(zq^{2M-2})^{-1})^{i_M,\ldots,i_{M+N}}_{i_M,\ldots,i_{M+N}}\big|_{i_{M}i_{M}}^{-1}
\\
&\times\big|(L^{\pm}(zq^{2M-2})^{-1})^{M+1,\ldots,M+N}_{M+1,\ldots,M+N}\big|_{M+1M+1}\ldots\big|(L^{\pm}(zq^{2M-2N})^{-1})^{M+N}_{{M+N}}\big|_{M+NM+N}.\\
\end{split}
\eeq
It is equal to $q^{-1}$-Berezinian of $\pi((L^{\pm}(zq^{2M-2N})^{-1})_{I^{c}})$.

The proof for $[1,M]\subseteq I $ is similar.
\epf

Taking $I=\emptyset$ in the above theorem, we immediately get that
\bco\label{ber inverse}
\beq
B_q(L^{\pm}(z))=B_{\qin} \left(\pi(L^{\pm}(zq^{2M-2N})^{-1}) \right).
\eeq
\eco

The following is an analogue of Schur's complement theorem (cf. \cite{Mo1}).
\bth\label{Sch com}
Write $L^{\pm}(z)$
in block matrix,
\[
L^{\pm}(z)=\left(\begin{array}{cc}
	L^{\pm}(z)_{11}&L^{\pm}(z)_{12}\\
	L^{\pm}(z)_{21}&L^{\pm}(z)_{22}\\
\end{array}\right)
\]
such that $L^{\pm}(z)_{11}$ and $L^{\pm}(z)_{22}$ are submatrices of size $k\times k$ and $(M+N-k)\times(M+N-k)$ respectively.
Then
\beq
B_q(L^{\pm}(z))=B_q(L^{\pm}(z)_{11}) B_{q} \left(L^{\pm}(zq^{2k})_{22}-L^{\pm}(zq^{2k})_{21}L^{\pm}(zq^{2k})_{11}^{-1}L^{\pm}(zq^{2k})_{12} \right)
\eeq
for $k<M$, and
\beq
B_q(L^{\pm}(z))=B_q(L^{\pm}(z)_{11}) B_{q}\left(L^{\pm}(zq^{4M-2k})_{22}-L^{\pm}(zq^{4M-2k})_{21}L^{\pm}(zq^{4M-2k})_{11}^{-1}L^{\pm}(zq^{4M-2k})_{12} \right)
\eeq
for $k\geq M$.
\eth
\bpf
Let
\[
X^{\pm}(z)=  L^{\pm}(z)^{-1}=\left(\begin{array}{cc}
	X^{\pm}(z)_{11}&X^{\pm}(z)_{12}\\
	X^{\pm}(z)_{21}&X^{\pm}(z)_{22}\\
\end{array}\right)
\]
It is well known that
\[
X^{\pm}(z)_{22}=(  L^{\pm}(z)_{22}-   L^{\pm}(z)_{21}   L^{\pm}(z)_{11}^{-1}   L^{\pm}(z)_{12})^{-1}.
\]
By Theorem \ref{Jac rat},
\beq
B_q(L^{\pm}(z))=B_q(L^{\pm}(z)_{11}) B_{\qin}  \left(\pi(X^{\pm}(z)_{22}) \right).
\eeq
It follows from Corollary \ref{ber inverse} that
\beq
B_{\qin}  \left(\pi(X^{\pm}(z)_{22}) \right) =B_{q} \\ \left(L^{\pm}(zq^{2k})_{22}-L^{\pm}(zq^{2k})_{21}L^{\pm}(zq^{2k})_{11}^{-1}L^{\pm}(zq^{2k})_{12} \right)
  \eeq
in the case $k< M$,
and
\beq
\bal
&B_{\qin}  \left(\pi(X^{\pm}(z)_{22}) \right) \\
&=B_{q}\left(L^{\pm}(zq^{4M-2k})_{22}-L^{\pm}(zq^{4M-2k})_{21}L^{\pm}(zq^{4M-2k})_{11}^{-1}L^{\pm}(zq^{4M-2k})_{12} \right)
  \eal
  \eeq
for $k\geq M$.

  This completes the proof.
\epf

Taking $k=M$ in Theorem \ref{Sch com} we have the following formula which is an analog of the definition of Berezinian for the classical super matrix. We note that it can also be obtained by Jacobi's ratio Theorem for the quantum affine algebra
$\mathrm U_{q}(\wh{\gl}_N)$.
\bco
If $k=M$ in Theorem \ref{Sch com}, then
\beq
B_q(L^{\pm}(z))={\det}_q(L^{\pm}(z)_{11}) {\det}_{q^{-1}} \left(L^{\pm}(zq^{2M})_{22}-L^{\pm}(zq^{2M})_{21}L^{\pm}(zq^{2M})_{11}^{-1}L^{\pm}(zq^{2M})_{12} \right)^{-1}.
\eeq
\eco

For any $k\geq 0$ we introduce the homomorphism
\beq
\phi_k: \U_{q}(\gl_{M|N})\rightarrow \mathrm \U_{q}(\gl_{k+M|N}),
\eeq
which takes $l^{\pm}_{ij}(z)$ to $l^{\pm}_{k+i,k+j}(z)$. Consider the composition
\beq
\psi_k=\om_{k+M|N}\circ\phi_k\circ\om_{M|N}.
\eeq

\ble\label{lemma syl}
For any $1\leq i,j\leq M+N$, we have
\beq
\psi_k(l^{\pm}_{ij}(z))=\left|\begin{array}{cccc}
	l^{\pm}_{11}(z)&\cdots&l^{\pm}_{1k}(z)&l^{\pm}_{1,k+j}(z)\\
	\vdots&&\vdots&\vdots\\
	l^{\pm}_{k1}(z)&\cdots&l^{\pm}_{kk}(z)&l^{\pm}_{k,k+j}(z)\\
	l^{\pm}_{k+i,1}(z)&\cdots&l^{\pm}_{k+i,k}(z)&\framebox{$l^{\pm}_{k+i,k+j}(z)$}\\
\end{array}\right|
\eeq
\ele

The proof of this lemma is almost the same as the case of Yangian (or super Yangian) (see Lemma 1.11.2 in \cite{Mo1}), so we omit it here.

The following Sylvester theorem for the quantum Berezinian follows from Lemma \ref{lemma syl} and Schur's complement theorem.
\bth\label{syl}
Denote the matrix of $\mathrm U_q(\wh{\gl}_{k+M|N})$ by $L^{\pm}(z)_{k+M|N}$, then
\[
\psi_k(B_q(L^{\pm}(z)))={\det}_q(L^{\pm}(zq^{-2k})^{1\ldots k}_{1\ldots k})^{-1}B_q(L^{\pm}(zq^{-2k})_{k+M|N}).
\]
\eth


The following is an analogue of the MacMahon Master theorem \cite{MRa} for the quantum affine superalgebra $\mathrm{U}_q(\widehat{\gl}_{M|N})$.
\bth
We have the identity
\beq
\sum_{r=0}^k(-1)^r\str\mathcal{S}^q_r\mathcal{A}^q_{\{r+1,\ldots,k\}}L^{\pm}_1(z)L^{\pm}_2(zq^2)\cdots L^{\pm}_k(zq^{2k-2})=0,
\eeq
\beq
\sum_{r=0}^k(-1)^r\str\mathcal{A}^q_r\mathcal{S}^q_{\{r+1,\ldots,k\}}L^{\pm}_1(z)L^{\pm}_2(zq^2)\cdots L^{\pm}_k(zq^{2k-2})=0,
\eeq
where $\mathcal{A}^q_{\{r+1,\ldots,k\}}$ and $\mathcal{S}^q_{\{r+1,\ldots,k\}}$ denote the antisymmetrizer and symmetrizer over the copies of $\End (\BC^{M|N})$ labeled by $r+1,\cdots,k$.
\eth
\bpf
We employ the quantum analog of the approach in \cite{MRa}.
It is sufficient to show that
\beql{mac}
\begin{split}
	&\str_{1,\cdots,k}\mathcal{S}^q_r\mathcal{A}^q_{\{r+1,\ldots,k\}}L^{\pm}_1(z)L^{\pm}_2(zq^2)\cdots L^{\pm}_k(zq^{2k-2})\\
	=&\str_{1,\cdots,k}\frac{[r]_q[k-r+1]_q}{[k]_q}\mathcal{S}^q_r\mathcal{A}^q_{\{r,\ldots,k\}}L^{\pm}_1(z)L^{\pm}_2(zq^2)\cdots L^{\pm}_k(zq^{2k-2})\\
	&+\str_{1,\cdots,k}\frac{[r+1]_q[k-r]_q}{[k]_q}\mathcal{S}^q_{r+1}\mathcal{A}^q_{\{r+1,\ldots,k\}}L^{\pm}_1(z)L^{\pm}_2(zq^2)\cdots L^{\pm}_k(zq^{2k-2})
\end{split}
\eeq

But due to \eqref{idem} the cyclic property, the RLL relation and the fact that $\mathcal{A}_r^q,\mathcal{H}_r^q $ are idempotents. It holds that
\[
\begin{split}
	&\str_{1,\cdots,k}\mathcal{S}^q_r\mathcal{A}^q_{\{r,\ldots,k\}}L^{\pm}_1(z)L^{\pm}_2(zq^2)\cdots L^{\pm}_k(zq^{2k-2})\\
	=&\str_{1,\cdots,k}\mathcal{S}^q_{r}\mathcal{A}^q_{\{r+1,\ldots,k\}}\frac{\hat{R}_{r,r+1}(q^{2r-2k})}{q^{k-r+1}-q^{r-k-1}}\mathcal{A}^q_{\{r+1,\ldots,k\}}L^{\pm}_1(z)L^{\pm}_2(zq^2)\cdots L^{\pm}_k(zq^{2k-2})\\
	=&\str_{1,\cdots,k}\mathcal{S}^q_{r}\frac{\hat{R}_{r,r+1}(q^{2r-2k})}{q^{k-r+1}-q^{r-k-1}}\mathcal{A}^q_{\{r+1,\ldots,k\}}L^{\pm}_1(z)L^{\pm}_2(zq^2)\cdots L^{\pm}_k(zq^{2k-2})\mathcal{A}^q_{\{r+1,\ldots,k\}}\\
	=&\str_{1,\cdots,k}\mathcal{S}^q_{r}\frac{\hat{R}_{r,r+1}(q^{2r-2k})}{q^{k-r+1}-q^{r-k-1}}\mathcal{A}^q_{\{r+1,\ldots,k\}}L^{\pm}_1(z)L^{\pm}_2(zq^2)\cdots L^{\pm}_k(zq^{2k-2})
\end{split}
\]
Similarly, we also have
\[
\begin{split}
	&\str_{1,\cdots,k}\mathcal{S}^q_{r+1}\mathcal{A}^q_{\{r+1,\ldots,k\}}L^{\pm}_1(z)L^{\pm}_2(zq^2)\cdots L^{\pm}_k(zq^{2k-2})\\
	=&\str_{1,\cdots,k}\mathcal{S}^q_{r}\frac{\hat{R}_{r,r+1}(q^{2r})}{q^{r+1}-q^{-r-1}}\mathcal{S}^q_{r}\mathcal{A}^q_{\{r+1,\ldots,k\}}L^{\pm}_1(z)L^{\pm}_2(zq^2)\cdots L^{\pm}_k(zq^{2k-2})\\
	=&\str_{1,\cdots,k}\frac{\hat{R}_{r,r+1}(q^{2r})}{q^{r+1}-q^{-r-1}}\mathcal{S}^q_{r}\mathcal{A}^q_{\{r+1,\ldots,k\}}L^{\pm}_1(z)L^{\pm}_2(zq^2)\cdots L^{\pm}_k(zq^{2k-2})
\end{split}
\]
Recall the equation $R-R'=(q-\qin)P$. Then we have that
\[
\frac{[r]_q}{[k]_q}\hat{R}_{r,r+1}(q^{2r-2k})+\frac{[k-r]_q}{[k]_q}\hat{R}_{r,r+1}(q^{2r})=q-\qin.
\]
These imply \eqref{mac}. Thus we have  the first equation. The second equation can be proved by the same arguments.
\epf

\bigskip
\centerline{\bf Acknowledgments}
\medskip

The work is supported in part by the National Natural Science Foundation of China (grant nos.
12171303 and 12001218), the Simons Foundation (grant no. MP-TSM-00002518),
and the Fundamental Research Funds for the Central Universities (grant nos. CCNU22QN002 and CCNU22JC001).		
\bigskip
\bigskip

\centerline{\bf Statement on Conflict of Interest}
\medskip
The authors have no conflict of interest to declare that are relevant to this article.
\bigskip
\bigskip

\centerline{\bf Statement on Data Availability }
\medskip
Data sharing is not applicable to this article as no new data were created or analyzed in this study.
\bigskip
\bigskip

\bibliographystyle{amsalpha}

\begin{thebibliography}{ABC}

\bibitem{BK} L.~Bagnoli and S.~Kozic, {\em A note on the quantum Berezinian for the double Yangian of the Lie superalgebra $\mathfrak {gl}_{m| n} $}, Algebr. Represent. Theory 28 (2025): 143-155.  

\bibitem{CH} H.~Chang and H.~Hu, {\em A note on the center of the super Yangian $Y_{M|N}(\mathfrak{s})$}, J. Algebra 633 (2023): 648-665.

\bibitem{Ch} I.~Cherednik,
{\em A new interpretation of Gelfand-Zetlin bases},
Duke Math. J.
54 (1987): 563-577.

\bibitem{GZ} M.~D.~Gould and Y.-Z. Zhang, {\em On super-RS algebra and Drinfeld realization of quantum affine superalgebras}, Lett. Math. Phys. 44 (1998): 291-308.

\bibitem{DF} J.~Ding and I.~B.~Frenkel, {\em Isomorphism of two realizations of quantum affine algebra $U_ q(\widehat{\mathfrak{gl}}(n))$}, Comm. Math. Phys.
156 (1993): 277-300.

\bibitem{FHS} H. Fan, B. Hou and K. Shi, {\em Drinfeld constructions of the quantum affine superalgebra
	$U_q (\widehat{\gl}(m|n))$}, J. Math. Phys, 38 (1) (1997): 411-433.

\bibitem{G} L.~Gow, {\em On the Yangian $Y(gl_{m|n})$ and its quantum Berezinian}, Czechoslovak J. Phys. 55 (2005): 1415-1420.

\bibitem{Jim}
{M.~Jimbo}, {\em A $q$-analogue of $U(\gl(N+1))$, Hecke algebra
	and the Yang--Baxter equation}, Lett. Math. Phys. {\bf 11} (1986):
247-252.

\bibitem{GTL} S.~Gautam and V.~Toledano~Laredo, {\em Yangians and quantum loop algebras}, Selecta Math. (N.S.) 19 (2) (2013): 271-336.

\bibitem{GR}
I.~Gelfand and V.~Retakh,
{\em Quasideterminants. {I}}, Selecta Math. (N.S.) 3 (1997): 517-546.

\bibitem{JLM}
N.~Jing, M.~Liu, and A.~Molev, {\em Eigenvalues of quantum Gelfand invariants}, J. Math. Phys. 65 (2024), 061703: (9pp).

\bibitem{KL}
D.~Krob and B.~Leclerc, {\em Minor identities for quasi-determinants and quantum determinants}, Comm. Math. Phys. 169 (1995): 1-23.

\bibitem{KwL} J.-H.~Kwon and S.-M.~Lee, {\it Super duality for quantum affine algebras of type $A$},
Int. Math. Res. Not. IMRN 2022 (2022): 18446-18525.

\bibitem{LRT} M.~de~Leeuw, V. Regelskis and A. Torrielli, {\it The quantum affine origin of the AdS/CFT secret symmetry}, J. Phys. A: Math. Theor. 45 (2012), 175202: (20pp).

\bibitem{LWZ} H.~Lin, Y.~Wang, and H.~Zhang, {\it From quantum loop superalgebras to super Yangians}, J. Algebra 650 (2024): 299-334.

\bibitem{L} K.~Lu,
{\it Gelfand-{Tsetlin} bases of representations for super {Yangian} and quantum affine superalgebra}, Lett. Math. Phys.
111 (2021), 145: (30pp).

\bibitem{LM} K.~Lu and E.~Mukhin, {\em Jacobi-Trudi identity and Drinfeld functor for super Yangian}, Int. Math. Res.
Not. IMRN 2021 (2021): 16751-16810.

\bibitem{Mo1}
{A.~Molev},
{\em Yangians for classical Lie algebras}.
Mathematical Surveys and Monographs
143. AMS, Providence, RI, 2007.

\bibitem{Mo2}
{A. Molev},
{\em Sugawara operators for classical Lie algebras}.
Mathematical Surveys and Monographs
229. AMS, Providence, RI, 2018.

\bibitem{MRa} A.~Molev and E.~Ragoucy, {\em  The MacMahon Master Theorem for right quantum superalgebras and higher Sugawara operators for $\widehat {\mathfrak {gl}}_{m|n} $}, Moscow Math. J. 14.1 (2014): 83-119.

\bibitem{MR} A.~Molev and V. Retakh, {\em Quasideterminants and Casimir elements for the general linear Lie superalgebra}, Int. Math. Res. Not. 2004 (2004): 611-619.

\bibitem{Na} M.~Nazarov, {\em Quantum Berezinian and the classical Capelli identity}, Lett. Math. Phys. 21 (1991): 123-131.

\bibitem{Na2} M.~Nazarov, {\em Yangian of the general linear Lie superalgebra}, SIGMA 16 (2020) 122: (24pp).

\bibitem{RS} N.~Yu.~Reshetikhin and M.~A.~Semenov-Tian-Shansky, {\em Central extensions of quantum current groups}, Lett. Math. Phys. 19 (1990): 133-142.

\bibitem{P} Y.-N.~Peng, {\em Parabolic presentations of the super Yangian $Y(gl_{M|N})$ associated with arbitrary 01-sequences}, Comm. Math. Phys. 346 (2016): 313-347.

\bibitem{T} A.~Tsymbaliuk, {\it Shuffle algebra realizations of type A super Yangians and quantum affine superalgebras for all Cartan data},
Lett. Math. Phys. 110 (2020): 2083-2111.

\bibitem{Y} H.~Yamane, {\em On defining relations of affine Lie superalgebras and affine quantized universal enveloping
	superalgebras}, Publ. Res. Inst. Math. Sci. 35 (1999): 321-390.

\bibitem{Z} H.~Zhang, {\em RTT realization of quantum affine superalgebras and tensor products},  Int. Math. Res. Not. IMRN 2016 (2016), No. 4: 1126-1157.

\bibitem{Zh} Y.-Z.~Zhang, {\em Super Yangian double and its central extension}, Phys. Lett., A 234 (1997): 20-26.


\end{thebibliography}

\end{document}